%
\input amstex
\documentstyle{amsppt}
\def\thedate{23 May 2011}
\TagsOnRight

\define\Muir{Mu}
\define\Barnard{Ba}
\define\Bourbaki{Bo}
\define\LakEke{EL}
\define\LakTho{LT}
\define\Macdonald{Mc}
\define\PohstZ{PZ}


\define\pref#1{{\rm(#1)}}
\define\Discr{\operatorname{Discr}}
\define\Vandermonde{\Delta}
\define\Ann{\operatorname{Ann}}
\define\Ker{\operatorname{Ker}}
\define\sign{\operatorname{sign}}
\define\Split{\operatorname{Split}}
\define\Sym{\goth S}

\def\sectno#1{\relax}

\topmatter
  \date \thedate \enddate
 \title On the invariants of the splitting algebra \endtitle
  \author Anders Thorup  \endauthor
 \affil Department of Mathematical Sciences, 
  University of Copenhagen\endaffil
 \address \rm Department of Mathematical Sciences, 
   Universitetsparken 5, DK--2100
     Copenhagen, Denmark \endaddress
 \email thorup\@math.ku.dk\endemail
 \abstract\nofrills{\bf Abstract.}\ \ 
 For a given monic polynomial $p(t)$ of degree $n$ over a commutative
ring $k$, the {\it splitting algebra\/} is the universal $k$-algebra
in which $p(t)$ has $n$ roots, or, more precisely, over which $p(t)$
factors,
 $$p(t)=(t-\xi_1)\cdots(t-\xi_n).
 $$
 The symmetric group $\Sym_r$ for $1\le r\le n$ acts on the splitting
algebra by permuting the first $r$ roots $\xi_1,\dots,\xi_r$. We give
a natural, simple condition on the polynomial $p(t)$ that holds if and
only if there are only trivial invariants under the actions. In
particular, if the condition on $p(t)$ holds then the elements
of $k$ are the only invariants under the action of $\Sym_n$.  

We show that for any $n\ge 2$ there is a polynomial $p(t)$ of degree
$n$ for which the splitting algebra contains a nontrivial element
invariant under $\Sym_n$. The examples violate an assertion by A. D.
Barnard from 1974. 
 \endabstract
 \endtopmatter

 \document

\subhead \sectno.1. Introduction \endsubhead
  Consider commutative algebras over a fixed commutative ring $k\ne
0$.  Fix a monic polynomial $p(t)$ of degree $n\ge 1$ with
coefficients in $k$:
 $$p(t)=a_0t^n+a_1t^{n-1}+\cdots+a_1t+a_n, \quad a_0=1.
\tag\sectno.1.1
 $$
 For $r=0,1,\dots,n$ let $\Split^{r}(p) =\Split^{r} (p/k)$ be the {\it
$r$'th splitting algebra\/} of $p(t)$, universal with respect to
factorizations,
 $$p(t)=(t-\xi_1)\cdots (t-\xi_r )\tilde p(t),\tag\sectno.1.2
 $$
 with $r$ factors $t-\xi_j$. In other words, such a factorization
exists in $\Split^{r}(p) [t]$, with elements $\xi_1,\dots,\xi_r $ in
$\Split^{r}(p) $ and a polynomial $\tilde p(t)\in \Split^{r}(p) [t]$,
and if $A$ is any $k$-algebra over which $p(t)$ factors,
 $$p(t)=(t-\alpha_1)\cdots (t-\alpha_r)q(t) ,\tag\sectno.1.3
 $$
 then there is a unique $k$-algebra homomorphism $\Split^{r}(p) \to A$
such that $\xi_j\mapsto \alpha_j$ for $j=1,\dots,r$, and,
consequently, $\tilde p(t)$ is mapped to $q(t)$. The ({\it
complete}) {\it splitting algebra\/} of $p(t)$ is obtained
when $r=n$; then $\tilde p(t)$ in (\sectno.1.2)
and $q(t)$ in (\sectno.1.3) are
equal to $1$.  

Clearly, the $r$'th splitting algebra is generated by the $r$
universal roots $\xi_1,\dots,\xi_r$ in \pref{\sectno.1.2}.  It follows
from the construction of $\Split^{r}(p)$ in Section 2 that the natural map
$\Split^{r}(p)\to \Split^{n}(p)$ is an injection, identifying
$\Split^{r}(p)$ with
the subalgebra $k[\xi_1,\dots,\xi_r]$ of
$\Split^{n}(p)=k[\xi_1,\dots,\xi_n]$.

Let $\Psi_{p}$ be the element of $k$ defined as the product $\Psi_{p}
=\prod_{i<j}(\xi_i+\xi_j)$. The product is a symmetric polynomial in
the roots $\xi_1,\dots,\xi_n$, and hence $\Psi _p$ a polynomial in the
coefficients $a_j$ of $p(t)$. In particular, $\Psi _p\in k$. There is
a simple determinantal formula for $\Psi_{p}$, see Definition\sectno.4.

It follows from the universal property that the symmetric group
$\Sym_n$ acts on the complete splitting algebra $\Split^{n}(p)$ by permuting
the roots $\xi_1,\dots,\xi_n$.  Obviously, the elements of the algebra
$\Split^{r}(p)=k[\xi_1,\dots,\xi_r]$ are invariant (or fixed) under
the action of the subgroup $\Sym_{n-r}'$ consisting of permutations in
$\Sym_n$ fixing the numbers $1,\dots,r$. The main result of this paper
is the following characterization, part of Theorem \sectno.7.

 \proclaim{Result} 
 The invariants are trivial: $k[\xi_1,\dots,\xi
_n]^{\Sym_{n-r}'}=k[\xi_1,\dots,\xi_r]$ for $0\le r\le n$, if and
only if $\Ann_k\Psi_{p}\cap \Ann_k 2=(0)$.
  \endproclaim

In particular, the equality $k[\xi _1,\dots,\xi_n]^{\goth S_n}=k$
holds if $\Ann_k\Psi _p\cap \Ann_k2=(0)$. Barnard \cite{\Barnard, p\.
289} asserted the equality for $n\ge 3$ without any condition on
$p(t)$.  We show by a counterexample that the general assertion is
not true.

Let $p^{(r)}(t)$ be the polynomial $\tilde p(t)$ in
\pref{\sectno.1.2}. Then the factorization has the form,
 $$p(t)=p_r(t)p^{(r)}(t),\quad\text{ where }
p_r(t):=(t-\xi _1)\dots (t-\xi _r).\tag\sectno.1.4
 $$
 Let $K_r=\operatorname{Fact}^r(p)$ be the $k$-subalgebra of
$\Split^{n}(p)$ generated by the elementary symmetric polynomials in
$\xi _1,\dots,\xi _r$, or, equivalently, by the coefficients of
$p_r(t)$. Then both polynomials $p^{(r)}(t)$ and $p_r(t)$ have
coefficients in $K_r$.  Under the condition $\Ann_k\Psi_{p}\cap \Ann_k
2=(0)$ on $p(t)$, we prove in Proposition \sectno.11 that the elements
of $K_r$ are the only invariants of the action of $\Sym_r$ on the
$r$'th splitting algebra, that is,
$k[\xi_1,\dots,\xi_r]^{\Sym_r}=K_r$.

\subhead \sectno.2. Construction\endsubhead
 A construction of the complete splitting algebra is given in
\cite{\Bourbaki, p.\thinspace IV.67, \S 5} and in \cite{\PohstZ,
p.~30}. We recall here the recursive construction of the intermediate
algebras $S_r:=\Split^{r}(p) $ for $r=0,\dots,n$:

 Obviously, $S_0=\Split^{0}(p)=k$.  For $r=1$, the equation
(\sectno.1.3) holds if and only if $\alpha_1\in A$ is a root of
$p(t)$. The universal algebra in which $p(t)$ has a root is obtained
by {\it adjoining formally a root\/} of $p(t)$:
 $$S_1=\Split^{1}(p):= k[x]/(p(x)), \qquad \xi_1:=(x\ \operatorname{mod}
p(t)).
 $$
 Assume that $S_r:=\Split^{{}r }(p/k)$ has been defined in general,
for $r<n$, with a factorization \pref{\sectno.1.2}.  Then, clearly, we
obtain $S_{r+1}$ by adjoining formally a root of $p^{(r)}(t)$, or,
equivalently, as the $r$'th splitting algebra of $p^{(1)}(t)$ over
$S_1$:
 $$\Split^{{}r +1}(p/k):= \Split^{1}(p^{(r)}/S_r )=\Split^{r}(p^{(1)}/S_1).
 $$

\proclaim{\sectno.3. Proposition}
 The monomials $\xi_1^{i_1}\xi_2^{i_2}\cdots \xi_{r}^{i_r}$ where
$0\le i_\nu\le n-\nu $ for $\nu =1,\dots, r$ form a $k$-basis for the
$r$'th splitting algebra $\Split^{r}(p)$.  In particular,
$\Split^{r}(p)$ is free of rank $n(n-1)\cdots(n-r+1)$ as a $k$-module,
and $\Split^{n}(p)$ is free of rank $n!$.
 \endproclaim\demo{Proof}
 The assertion follows by induction on $r$ from the recursive
definition of $\Split^{r}(p)$.
\enddemo

\remark{Note}
 It is an easy consequence of the  Proposition that the roots
$\xi_1,\dots,\xi_n$ are $n$ different elements in $\Split^{n}(p)$ except
when $n=2$, $a_1=0$, and $2=0$ in $k$.
 \endremark

 \subhead \sectno.4. Definition\endsubhead
 Consider in $\Split^{n}(p)=k[\xi_1,\dots,\xi_n]$ the Vandermonde
determinant $\Delta_{p}$,
 $$\Delta_{p} =\prod_{i<j}(\xi _i-\xi _j) =\sum_{\sigma\in
\Sym_n}(\operatorname{sign} \sigma )\, \sigma\bigl(\xi _1^{n-1}\xi
_2^{n-2}\cdots \xi _{n-1}\bigr),
 $$
 and the two elements $\Psi_{p}$ and $\Discr_{p}$,
 $$\Psi_p= \prod_{i<j} (\xi_i+\xi_j)
 \quad\text{ and }\Discr_p=\prod_{i<j} (\xi_i-\xi_j)^2.
 $$
 The element $\Discr_{p}=\Delta_{p}^2$ 
is of course the {\it discriminant\/} of $p(t)$.

The elements $\Psi_{p}$ and $\Discr_{p}$ are symmetric polynomials in the
roots $\xi_1, \dots , \xi_n$, and consequently can be expressed
as polynomials in the coefficients of $p(t)$. In particular, the
elements $\Psi_{p}$ and $\Discr_{p}$ belong to $k$.

It is well known that $\Psi_{p}$, as a polynomial in the $\xi_i$, is a
Schur polynomial, see \cite{\Macdonald, Example 7, p.~46} or
\cite{\Muir, formula 339 p.~334}. As a polynomial in the elementary
symmetric polynomials $e_j$ it is the determinant $\Psi=\det
(e_{2i-j})$, see \cite{\Macdonald, Formula (3.5), p.~41}. In terms of
the coefficients $a_i=(-1)^ie_i$ of $p(t)$ we obtain the expression
$\Psi_{p} =(-1)^{n(n-1)/2}\det(a_{2i-j})$,
 $$\Psi_{p}= (-1)^{n(n-1)/2} \left|\matrix
a_0&0&0&0&\cdots&0\\ a_2&a_1&a_0&0&\cdots &0\\
a_4&a_3&a_2&a_1&\cdots&0\\
\vdots&\vdots&\vdots&\vdots&\ddots&\vdots\\
a_{2n-2}&a_{2n-3}&a_{2n-4}&a_{2n-3}&\cdots&a_{n-1}
 \endmatrix\right|.
 $$
 In particular, in low degrees: For $n=1$: $\Psi=1$, for $n=2$:
$\Psi=-a_1$, for $n=3$: $\Psi =a_3-a_1a_2$, and for $n=4$: $\Psi
=a_1a_2a_3 -a_1^{\,2} a_4 -a_3^{\,2}$.

\proclaim{\sectno.5. Lemma} Assume for $n\ge 2$ that  
$F\in k[\xi_1]$ is $\Sym_{2}$-invariant.  Then:

{\rm(1)} If $n\ge 3$, then $F\in k$. 

{\rm(2)} If $n=2$, then $F=b\xi _1+c$, with $b,c\in k$ and $2b=\Psi
_pb=0$.
 \endproclaim
 \demo{Proof}
 Let $\tau$ be the non-trivial permutation in $\Sym_2$, acting on the
second splitting algebra $k[\xi _1,\xi _2]$ by interchanging $\xi _1$
and $\xi_2$. By assumption, $F\in k[\xi_1]$ and $\tau F=F$.  Write $F$
in the form $F=Q(\xi_1)$, where $Q\in k[x]$ is a polynomial of degree
less than $n$, say $Q=bx^{n-1}+cx^{n-2}+\cdots\,$. Then
$Q(\xi_2)=Q(\xi_{1})$.  So the polynomial $Q(x)-Q(\xi_1)$ in
$k[\xi_1][x]$ has $\xi_{2}$ as a root.  Hence $Q(x)-Q(\xi_1)$ is a
multiple of $p^{(1)}(x)$.  As $p^{(1)}(x)$ is monic of degree $n-1$ it
follows by comparing the degrees and the leading coefficients that
 $$Q(x)-Q(\xi_1)=b p^{(1)}(x).\tag\sectno.5.1
 $$
 After multiplication by $x-\xi_1$, we obtain
the following equation in $k[\xi_1][x]$:
 $$\bigl(Q(x)-Q(\xi_1)\bigr)(x-\xi_1)=bp(x).\tag\sectno.5.2
 $$
 Compare the coefficients of $x^{n-1}$ in the equation.  If $n\ge 3$
we obtain the equation $-b\xi _1+c=ba_1$. In particular, $b\xi _1\in
k$, and hence $b=0$.  From \pref{\sectno.5.1} we conclude that the
polynomial $Q(x)$ is a constant.  Hence $F=Q(\xi_1)$ belongs to $k$.
Thus Part (1) has been proved.

The case $n=2$ is easily treated directly.  Alternatively we may use
\pref{\sectno.5.2}. Equating the coefficients of $x$ gives the equation
$-2b\xi _1=ba_1$. Here $a_1=-(\xi_1+\xi_2)=-\Psi _p$, and hence $2b\xi
_1=\Psi _pb$.  As $\xi _1,1$ are linearly independent over $k$, it
follows that $2b=\Psi _pb=0$. Thus Part (2) has been proved. 
 \enddemo

\proclaim{\sectno.6. Proposition}
 Assume for $n>r\ge 1$ that $F\in k[\xi _1,\dots,\xi _r]$ is
$\Sym_{r+1}$-invariant.  Then: 

{\rm(1)} If $r\le n{-}2$ then $F\in k$.

{\rm(2)} If $r=n-1$ then $2F\in k$ and $\Psi_{p}F\in k$. 
 \endproclaim
 \demo{Proof}
 Set $S_j:=k[\xi _1,\dots,\xi _j]$. Then $S_j=S_{j-1}[\xi _j]$ is the
first splitting algebra of $p^{(j-1)}(t)$ over $S_{j-1}$. The
degree of $p^{(j-1)}$ is $n-j+1$, and hence at least $3$ if $j\le n-2$. 

Therefore, under the assumptions in Part (1), it follows by repeated
application Lemma \sectno.5\,(1) that $F\in S_{j-1}$ for $j=r,\dots,1$.
With $j=1$ it follows that $F\in k$. 

For Part (2), note first that the assertion for $n=2$ follows from
Lemma \sectno.5\,(2).  Proceed by induction on $n\ge 3$. Note that
$\Split^n(p)$ is the complete splitting algebra of $p^{(1)}(t)$ over
$k[\xi_1]$. Clearly 
 $$\Psi_{p}=\Phi\Psi_{p^{(1)}}  \quad\text{where } 
\Phi:=\prod\nolimits_{1<j\le n}(\xi_1+\xi_j).
 $$
 Moreover, $\Phi\in k[\xi_1]$; in fact $\Phi=(-1)^{n-1}p^{(1)}(-\xi_1)$.
By induction, $2F$ and $\Psi_{p^{(1)}}F$ belong to $k[\xi_1]$.  So both
products $2F$ and $\Psi_{p}F$ belong to $k[\xi_1 ]$. As both products
are $\Sym_n$-invariant, it follows from Lemma 6\,(1) that they belong
to $k$.
 \enddemo

\proclaim{\sectno.7. Theorem}
 Let $S=\Split^{n}(p/k)=k[\xi_1,\dots,\xi_n]$ be the complete splitting
algebra of $p(t)$. Assume that $n\ge 2$. Then the following
conditions on $p(t)$ are equivalent:
 \roster
 \item"(i)" $\Ann_k \Discr_{p}\cap \Ann_k 2=(0)$. 
 \item"(ii)" $\Ann_k \Psi_{p}\cap \Ann_k 2=(0)$.
  \item"(iii)" $S^{\Sym_2'} =k[\xi_1,\dots,\xi_{n-2}]$.
 \item"(iv)" $S^{\Sym_{n-r }'} =k[\xi_{1}, \dots , \xi_{r }]$ for
$r=0,1, \dots , n-2$, where $\Sym_{n-r}'$ denotes the subgroup
of permutations in $\Sym_n$ fixing the numbers $1,\dots,r$. 
 \endroster
 \endproclaim
 \demo{Proof}
 For an element $F\in S$ and a subset $V\subseteq S$ denote by $F|V$
the restriction to $V$ of multiplication by $F$.  Consider with
$I:=\Ann_S2$ the following three conditions:
 $$\def\quad{\hskip0,8em }%
\text{(i*)}\ \Ker \Discr_p|I=(0),\qquad
  \text{(ii*)}\ \Ker \Psi_p|I=(0), \qquad
  \text{(iii*)}\ \Ker (\xi_{n-1}+\xi_n)|I=(0).
 $$
 The algebra $S$ is free as a $k$-module and the elements $\Psi_p$
and $\Discr_p$ belong $k$.  Hence (i) is equivalent to (i*) and
(ii) is equivalent to (ii*).  Set $S_{r-2}=k[\xi _1,\dots,\xi
_{r-2}]$. Then $S$ is the complete
splitting algebra over $k[\xi _1,\dots,\xi _{n-1}]$ of the degree
$2$-polynomial $p^{(n-2)}$.  By Lemma \sectno.5\,(2), (iii) holds if and
only if $\Ann_{S_{n-2}}2\cap \Ann_{S_{n-2}}(\xi _{n-1}+\xi _n)=(0)$.
Again, as $S$ is free over $S_{n-2}$, the latter condition holds if
and only if (iii*) holds.

Since $I=\Ann_S2$ it follows that $\Psi_p|I=\Vandermonde_p|I$.  Hence
$\Discr_p|I=(\Psi_p|I)^2$.  Consequently $\Psi_p |I$ is injective if and
only if $\Discr_p|I$ is injective.  Hence (i*)$\iff$(ii*). 

Again, $\Psi_p|I$ is the product of the factors $(\xi _i+\xi _j)|I$
for $i<j$. Hence, if $\Psi|I$ is injective, then the
factor $(\xi _{n-1}{+}\xi _n)|I$ is injective. Assume conversely that
the factor $(\xi _{n-1}{+}\xi _n)|I$ is injective.  Then, since the
group $\Sym_n$ acts by automorphisms of $S$, every factor
$(\xi_i{+}\xi_j)|I$ (with $i<j$) is injective, and hence the product
$\Psi|I$ is injective. Hence (ii*)$\iff$(iii*). 

Obviously (iii) is part of the conditions in (iv), and hence
(iv)$\Rightarrow $(iii). Clearly, to finish the proof it suffices to
show that (ii) implies that the equality in (iv) for $r=0$ holds.  So
assume that $F\in S$ is $\Sym_n$-invariant. Then, by Proposition
\sectno.6 we have the relations $2F=0$ and $\Psi_pF=0$. The relations
mean that if $F$ is expanded in terms of the basis in Proposition
\sectno.3 then all
coefficients to base elements different from $1$ are annihilated by
$2$ and by $\Psi_p$.  Therefore, if (ii) holds then all these
coefficients vanish, that is, $F\in k$.
 \enddemo

 \proclaim{\sectno.8. Corollary}
 If any of the three elements: $2$, or $\Discr_{p}$, or $\Psi_{p}$, is
a non-zero divisor in $k$, then the conditions of the Theorem hold. In
particular, then
 $$k[\xi_1,\dots, \xi_n]^{\Sym_n}=k.\tag\sectno.8.1
 $$
 \endproclaim\demo{Proof}
 Clearly, under the given assumptions either (i) or (ii) of the
Theorem holds.  Hence (iv) holds, and in particular \pref{\sectno.8.1}
which is the special case $r=0$ of (iv) holds.
  \enddemo

\subhead\sectno.9. Notes\endsubhead
 The results in the Corollary for the elements $2$ and $\Discr_{p}$
were proved by Pohst and Zassenhaus \cite{\PohstZ, (2.18d),
p.\thinspace46 and (3.6), p.\thinspace49} and, with a different proof,
by Laksov and Ekedahl \cite{\LakEke, Theorem 5.1 and Remark 5.3,
p.~13--14}. Pohst and Zassenhaus also proved the assertion in
Proposition \sectno.6\,(2), but with $\Psi _{p}$ replaced by
$\Discr_{p}$.

 It was asserted by Barnard \cite{\Barnard, Proposition 4, p\. 289}
that the equality \pref{\sectno.8.1} holds for all $n>2$.  However, a
simple counterexample shows that the assertion cannot hold in the
stated generality: Consider the splitting algebra of the polynomial
$p(t)=t^n$ ($n\ge 2$). It is easy to see, by induction on
$r=1,\dots,n-2$, that $\xi _1^{n-1}\cdots \xi _{r}^{n-r}\xi _r=0$.
 It follows easily for $\sigma\in \Sym _n$ that 
 $$\sigma \bigl(\xi_1^{n-1}\xi _2^{n-2}\cdots \xi _{n-1} \bigr)
= (\sign \sigma )\xi _1^{n-1}\xi
_2^{n-2}\cdots \xi _{n-1}.
 $$
  Hence for $z\in
\Ann_k2$ the element $z\xi _1^{n-1}\xi _2^{n-2}\cdots \xi _{n-1}$ is
invariant, and it is non-trivial if $z\ne0$. A second family of
non-trivial invariants is given in Example \sectno.10.

It is an open question, at least to the knowledge of the author,
whether the equality \pref{\sectno.8.1} implies the stronger
conditions in Theorem \sectno.7.

\subhead\sectno.10. Example\endsubhead A natural idea to construct
invariants in $\Split^{n}(p)$ ($n\ge 2$) is to write the Vandermonde
determinant $\Delta _p$ as the difference,
 $$\Delta_p=\Delta^+-\Delta^-, \qquad 
\Delta^+=\sum_{\sigma\in \goth A_n}
\sigma\bigl(\xi _1^{n-1}\xi _2^{n-2}\cdots \xi _{n-1}\bigr),
 $$
 where the sum is over all even permutations. Then $\Delta^+ $ and
$\Delta^-$ are invariant under even permutations and interchanged by
odd permutations.  In particular, if $z\in k$ then $z\Delta^+$ is
invariant under $\Sym_n$ if and only if $z$ is in the kernel of
multiplication by $\Delta_p$ as a map $k\to S$.  But naturally, even if
$z\ne 0$ and $z\Delta_p=0$, it may happen that $z\Delta ^+=0$. 

In particular, assume that $z\in \Ann_k2$. Then, as noted above,
$z\Delta_{p} =z\Psi_{p}$, and hence $z\Delta_{p}^+$ is invariant if
and only if $z\in \Ann_k\Psi_p$. For $n=2$ or $n=3$ it is easy to see
that the invariants of the complete splitting algebra are the elements
$c+z\Delta_{p} ^+$ for $c\in k$ and $z\in \Ann_k2\cap \Ann
_k\Psi_{p}$.

\proclaim{\sectno.11. Proposition}
 Fix $r$ with $1\le r\le n$ and let $K=\operatorname{Fact}^r(p)$ be
the $k$-subalgebra of $\Split^{n}(p)$ generated by the elementary
symmetric polynomials in the first $r$ roots $\xi _1,\dots,\xi_r$, or,
equivalently, by the coefficients of the polynomial
$p_r(t):=(t-\xi_1)\cdots(t-\xi_r)$.  Then, in $K[t]$ we have the
factorization,
 $$p(t)=p_r(t)p^{(r)}(t),\tag\sectno.11.1
 $$
   of $p(t)$ into two monic factors, the first of degree $r$ and the
second of degree $n-r$.  Moreover, the $k$-algebra $K$ is universal
with respect to this property.  Furthermore, the algebra
$\Split^{r}(p)=k[\xi_1,\dots,\xi_r]$ is the complete splitting algebra of
the degree-$r$ polynomial $p_r(t)$ over $K$:
 $$\Split^{r}(p/k)=\Split^{r}(p_r/K).\tag\sectno.11.2
 $$
 Finally, if the equivalent conditions of Theorem\/ {\rm \sectno.7} hold for
$p(t)$, then
 $$k[\xi_1,\dots,\xi_r]^{\Sym_r}=K.\tag\sectno.11.3
 $$
 \endproclaim
 \demo{Proof}
 The equation \pref{\sectno.11.1} is simply \pref{\sectno.1.2}.
The polynomial $p_r(t)$ has, by construction, coefficients in $K$.
Therefore, by \pref{\sectno.11.1}, so has $p^{(r)}(t)$.

Obviously, the polynomial $p_r(t)\in K[t]$ splits completely over
$S_r:=k[\xi_1,\dots,\xi_r]$. To prove \pref{\sectno.11.2}, we verify
that the splitting is universal. So assume that $\varphi_0\:K\to
A$ is an algebra such that $p_r(t)$ factors completely over $A$ with
$r$ factors $t-\alpha_j$.  Then we obtain in $A[t]$ the
factorization \pref{\sectno,1.3} where $q(t)$ is the image in $A[t]$ of
$p^{(r)}(t)$. So, by the universal property of $S_r=\Split^{r}(p/k)$,
there is a unique $k$-algebra homomorphism $\varphi \:S_r\to A$ such
that $\varphi (\xi_j)=\alpha_j$ for $j=1,\dots,r$. It remains to
prove that $\varphi $ is a map $K$-algebras, that is, $\varphi $ is
equal to $\varphi_0$ on the $k$-subalgebra $K$ of $S_r$.  The
equality results from the fact that under both maps the coefficients
of $p_r(t)$ are mapped to the signed elementary symmetric polynomials
of the $\alpha_j$. So the two maps agree on the coefficients of
$p_r(t)$, and since $K$ is generated as a $k$-algebra by these
coefficients, the two maps agree on $K$.

The universal property of $K$ with respect to factorizations
$p(t)=\tilde q(t)q(t)$ with two factors of degrees $r$ and $n-r$ is
proved similarly: Assume that such a factorization exists over a
$k$-algebra $A$.  Since $K$ is generated by the coefficients of
$p_r(t)$ there is at most one $k$-algebra homomorphism $K\to A$
under which $p_r(t)$ is mapped to $\tilde  q(t)$.  To prove the
existence, consider the complete splitting algebra $T_r$ of $\tilde 
q(t)$ over $A$. Then there is a $k$-algebra homomorphism $S_r\to T_r$
such that the $\xi_j$ are mapped to the roots of $\tilde  q(t)$. In
particular, $p_r(t)$ is mapped to $\tilde  q(t)$. As $K$ is
generated by the coefficients of $p_r(t)$, and the coefficients of
$\tilde  q(t)$ belong to $A$, we obtain the map $K\to A$ as the
restriction of the map $S_r\to T_r$.

To prove the final assertion, consider the equivalent conditions of
Theorem \sectno.7. Assume that they hold for $p(t)$. Note that the
factors defining the product $\Psi _{p_r}$ are also factors of $\Psi
_p$. Therefore,
 $$\Ann_K2\cap \Ann_K\Psi _{p_r}\subseteq \Ann_{S_r}2\cap \Ann_{S_r}\Psi _{p}.
 $$
 Since ${S_r}$ is free over $k$, condition (ii) for $p(t)$ implies
that the right hand intersection is trivial.  Therefore the left hand
intersection is is trivial, that is, condition (ii) holds for $p_r(t)$
in $K[t]$. Moreover, $k[\xi_1,\dots,\xi _r]$ is the complete splitting
algebra of $p_r(t)$ over $K$ by \pref{\sectno.11.2}. Therefore, by the
Theorem, condition (iv) holds for $p_r(t)$; in particular
(\sectno.11.3) holds.
 \enddemo

\subhead\sectno.12. Note\endsubhead
 The algebra $\Split^{r}(p)$ is free over $k$ of rank $n(n-1)\cdots
(n-r+1)$ by construction, and it is free over
$K=\operatorname{Fact}^r(p)$ of rank $r!$ by \pref{\sectno.11.2}.  We
showed in the paper with D. Laksov \cite{\LakTho} that $K=\operatorname{Fact}^r(p)$
is in fact $k$-free of rank $\binom nr$, generated by suitable Schur
polynomials in $\xi _1,\dots,\xi _r$.  The paper describes in addition
the connection between splitting algebras and intersection rings of
Grassmannians (Schubert Calculus).

\enddocument